\documentclass[12pt,reqno]{amsart}

\usepackage{mathrsfs}
\usepackage{amsmath}
\usepackage{latexsym}
\usepackage{colortbl}
\usepackage{amssymb}
\usepackage{amsfonts}
\usepackage{graphics, hyperref}
\include{PDF}

\renewcommand{\proof}{\par\noindent{\it Proof.\ \ }}
\def\qed{\ifmmode\square\else\nolinebreak\hfill
$\Box$\fi\par\vskip12pt}

\def\ov{\overline} 
\def\l{\langle} \def\r{\rangle} 
\def\div{\,\big|\,} 

\def\FF{\mathbb F}

\def\ZZ{{\rm C}}

\def\S{{\rm S}}

\def\C{{\bf C}} \def\N{{\bf N}}
\def\Z{{\bf Z}}

\def\mod{{\sf mod~}}

\def\Aut{{\sf Aut}} \def\Inn{{\sf Inn}} \def\Out{{\sf Out}}

\def\char{\,{\sf char}\,}

\def\a{\alpha} \def\b{\beta} \def\d{\delta} \def\s{\sigma}
 \def\o{\omega} 

\def\Hom{{\rm Hom}}
\def\ord{{\rm ord}}
\def\Irr{{\rm Irr}}

\def\GammaL{{\rm \Gamma L}}

\def\Sym{{\rm Sym}}

\def\GL{{\rm GL}}  
\def\AGL{{\rm AGL}}

\newtheorem{theorem}{Theorem}[section]%
\newtheorem{lemma}[theorem]{Lemma}%
\newtheorem{corollary}[theorem]{Corollary}%
\newtheorem{proposition}[theorem]{Proposition}%
\newtheorem{example}[theorem]{Example}%

\begin{document}
\title[Automorphism groups]
{On the automorphism groups of Frobenius Groups}

\thanks{This work was supported by NSF of Yunnan Province(Grant No. 2017FD071).}

\author{Lei Wang}
\address{Lei Wang School of Mathematics and Statistics\\
Yunnan University\\
Kunming, Yunnan 650091\\
P. R. China}

\email{wanglei@ynu.edu.cn}

\date{\today}



\maketitle

\begin{abstract}
This is one of a series papers which aim towards to
solve the problem of determining automorphism groups of Frobenius groups.
This one solves the problem in the case where the
Frobenius kernels are elementary abelian and
Frobenius complements are cyclic.

\end{abstract}

\qquad {\textsc k}{\scriptsize \textsc {eywords.}} {\footnotesize
Frobenius group, automorphism group}

\section{Introduction}

A Frobenius group $G$ is a semidirect product of a normal subgroup $V$ by a subgroup $H$
such that none of the non-identity elements of $H$ centralizes a non-identity element of $V$,
where $V$ is called the {\it Frobenius kernel} and $H$ is called a {\it Frobenius complement} of $G$.
Furthermore, by the well-known result of Thompson,
the Frobenius kernel is a nilpotent group, and
by Burnside's result,
each Sylow subgroup of a Frobenius complement is a metacyclic group with restricted properties.

Frobenius groups form an important class of groups, and
have been extensively studied in the literature,
refer to \cite{DM-book,Gorenstein,B,Sta,Thomp1,Thomp2,Thomp3,Zas}.
A natural problem arises:

\vskip0.1in
{\noindent\bf Problem A.} Determine automorphism groups of Frobenius groups.
\vskip0.1in

The problem is unsolved yet.
This paper solves it in the case where the Frobenius kernels
are elementary abelian, and Frobenius complements are cyclic.

Frobenius groups have played an important role
not only in group theory, but also in various applications,
refer to \cite[p.320-324]{Jac} for the applications
in algebraic structures, and
see \cite{Corr,JinWei,Fang,Gur,L-P,Li-Wang,Song,S}
for  the applications in algebraic graph theory.
In these applications, determining automorphism groups
of certain Frobenius groups is a crucial step.
This is actually one of our principle motivations
for the work of this paper.

In order to state our results, we need to introduce some notation.
For a finite group $G$ and a field $\FF$, we denote by
$\Irr(\FF G)$ a complete set of representatives
for the isomorphism classes of irreducible $\FF G$-modules.
But we mean by $V\in\Irr(\FF G)$ that $V$ is an irreducible
$\FF G$-module for convenience.
For a positive integer $e$ and an
element $V\in\Irr(\FF G)$,
we abuse notation and denote by $V^e$ a homogeneous $\FF G$-module,
which is a direct sum of $e$ copies of $V$.
Let $H\leqslant G$, and let $V$ be an $\FF G$-module.
We denote by $V_H$ the $\FF H$-module
obtained by restricting the action of $\FF G$ on $V$ to $\FF H$.
For positive integers $a,n$,
we call $m$ the order of $a$ modulo $n$ if
$n$ divides $a^m-1$ but
$n$ does not divide $a^i-1$ for $i<m$,
and denote $m$ by $\ord_n(a)$.

\begin{theorem}\label{WLLLL}
Let $G=V{:}H=\ZZ_p^d{:}\ZZ_n$ be a Frobenius group,
where $p$ is a prime, and $d,n$ are positive integers.
Then

\[\Aut(G)=V.(((\GL(e_1,p^f)^\ell\times\dots\times
\GL(e_s,p^f)^\ell).\ZZ_f).L),\]
where $\ZZ_f.L\leqslant\Aut(H)$, $f=\ord_n(p)$, $|L|=\ell$,
and $(e_1+\dots+e_s)f\ell=d$.
\end{theorem}
\vskip0.1in

\noindent{\bf Remarks on Theorem~\ref{WLLLL}.}
\begin{itemize}
\item[(a)]
The automorphism group
$\Aut(G)=V.\N_{\Aut(V)}(H)$
by Lemma~\ref{Aut(G)}.

\vskip0.1in
\item[(b)]
Let $M{:}=\N_{\Aut(V)}(H)$. 
View $V$ as an $\FF_pM$-module.
By Lemma~\ref{irrede},
\[\mbox{$V_M=V_1\oplus V_2\oplus\cdots\oplus V_s$,
where $V_i\in\Irr(\FF_pM)$
and $V_i\not\cong V_j$ $(i\not=j)$}.\]

\vskip0.1in
\item[(c)]
By Lemma~\ref{FM}, $V_i$ can be decomposed as
\[\mbox{$(V_i)_H=V_{i1}^{e_i}\oplus
V_{i2}^{e_i}\oplus\cdots\oplus V_{i\ell}^{e_i}$,}\]
where $V_{ij}\in\Irr(\FF_pH)$, $\dim_{\FF_p}V_{ij}=f$
and $V_{ij}\not\cong V_{ik}$ $(j\not=k)$.

\vskip0.1in
\item[(d)]
The automorphism group $\ZZ_f$ is called the field automorphism
group of $H$, which is induced  by the Frobenius automorphism
of $\FF_{p^f}$.

\end{itemize}

\vskip0.1in
In what follows, we will adopt the convention:
the Frobenius kernel $V$ 
is viewed as the vector space $\FF_p^d$
if $V$ is considered as an $\FF_pM$-module,
and an elementary abelian group $\ZZ_p^d$ otherwise.
This conventional device will be extremely useful in passing between the group-theoretic and representation theory points of view.
\vskip0.1in

In subsequent work, we apply Theorem~\ref{WLLLL}
to characterise a class of finite groups,
and their Cayley graphs.

\section{Preliminary}\label{Cayley-sec}

We first define some notation.
For a group $G$, denote by $\Z(G)$ the center of $G$.
For a group $T$ and
an integer $\ell$, by $T^\ell$
we mean the direct product of $\ell$ copies of $T$.
For a positive integer $n$ and a prime $p$, denote
by $\ZZ_n$ and $\ZZ_p^n$ a cyclic group of order $n$
and an elementary abelian group of order $p^n$.
Given two groups $N$ and $K$, denote by $N\times K$
the direct product of $N$ and $K$, by $N.K$ an extension of $N$ by $K$,
and if such an extension is split, then we write $N{:}K$
instead of $N.K$.

\begin{lemma}\label{Inn}
Let $G$ be a finite group with $\Z(G)=1$. Then $\C_{\Aut(G)}(\Inn(G))=1$.
\end{lemma}
\proof By our assumption, $G\cong\Inn(G)$.
Let $\phi$ be an isomorphism from $G$ to $\Inn(G)$.
Then each element $\phi(g)$ of $\Inn(G)$ acts on $G$ by
$x^{\phi(g)}=g^{-1}x g$, where $g,x\in G$.

Let $L:=\C_{\Aut(G)}(\Inn(G))$.
Then for $\s\in L$, we have
\[\mbox{$\phi(g)^{\s}=\s^{-1}\phi(g)\s=\phi(g)$,}\]
and hence
\[\mbox{$z^{\s^{-1}\phi(g)\s}=z^{\phi(g)}$ for any $z\in G$.}\]
Since $z^{\s^{-1}\phi(g)\s}=(g^{\s})^{-1}zg^{\s}$,
it follows that $(g^{\s})^{-1}zg^{\s}=g^{-1}zg$ by the
previous equation, and so $z^{g^{\s}g^{-1}}=z$. Thus
$g^{\s}g^{-1}$ centralises $G$.
So $g^{\s}g^{-1}\in\Z(G)$, forcing $g^\s=g$.
Since $g$ is arbitrary,
it implies that $\s$ fixes $G$ pointwise, and so $\s=1$.
Thus $L=1$. This completes the proof.
\qed

Let $G=N{:}H$ be a finite group,
where $\C_H(N)=1$, $\Z(G)=1$
and $\gcd(|N|,|H|)=1$.
By Lemma~\ref{Inn},
we will identify $G$ with $\Inn(G)$ a normal subgroup of $\Aut(G)$.

\begin{lemma}\label{Centr}
$\C_{\Aut(G)}(N)=\Z(N)$.

\end{lemma}
\proof Let $A:=\Aut(G)$.
By the convention made above, we have
\[N\,\char\, N{:}H=G\lhd A,\]
and hence $N\lhd A$.
Let $C:=\C_A(N)$.
Then $\Z(N)\leqslant C$.
We will show that in fact
\[\mbox{$C=\Z(N)$.}\]

Let $\ov G=G/\Z(N)$ and $\ov C=C/\Z(N)$.
Since $N \char G\lhd A$, we have $\C_G(N)\unlhd A$.
By our assumption, $\C_G(N)\leqslant N$,
and so $\C_G(N)=\Z(N)$.
It is clear that $G C/\Z(N)=\ov G\times \ov C$.
We pick any element
\[\theta\in C=\C_A(N).\]
Notice that for $g\in G$, we have
$[\theta,g^{-1}]\in\Z(N)$.
It follows that $g^{\theta}\in\Z(N)g$,
that is,
\[\mbox{$g^{\theta}=x_g g$ where $x_g\in\Z(N)$.}\]
If $\Z(N)=1$, then $g^{\theta}=g$.
Since $g$ is arbitrary, we have $\theta=1$, and so
$C=1$.  Therefore, we may assume that $\Z(N)\not=1$.

Let $R:=\Z(N){:}H$. By the previous paragraph,
we deduce that $R^\theta=R$,
and thus $H^\theta\leqslant R$.
Since $\gcd(|N|,|H|)=1$, it follows that
$H$ and $H^\theta$ are two Hall subgroups of $R$.
By Schur-Zassenhaus's Theorem, there exists some
$y\in\Z(N)$ such that $H^{\theta y^{-1}}=H$.
Note that $\Z(N)\cap H=1$.
It is easily shown that
$\theta y^{-1}$ centralises $H$.
Clearly, $\theta y^{-1}$ centralises $N$.
It follows that $\theta y^{-1}$ centralises $G$.
Thus $\theta=y$. Since $\theta$ is arbitrary,
it implies that $C=\Z(N)$, completing the proof.
\qed

In what follows, we continue to use the notation above.
By Lemma~\ref{Centr},
\[G/\Z(N)=\ov G\unlhd\ov A:=A/\Z(N)\lesssim\Aut(N).\]
Next identify $\ov G$ and $\ov A$ with the
subgroups of $\Aut(N)$.
Since $\gcd(|N|,|H|)=1$,
it follows that $\Inn(N)$ \char\, $\ov G$,
and so $\Inn(N)\unlhd\ov A$.
Let
\[\mbox{$\ov{\ov A}:=\ov A/\Inn(N)$
and  $\ov{\ov H}:=\ov G/\Inn(N)$}.\]
Then $\ov{\ov H}\lhd\ov{\ov A}\leqslant\Out(N)$,
and thus $\ov{\ov A}\leqslant\N_{\Out(N)}(\ov{\ov H})$.
Furthermore, we have

\begin{lemma}\label{Aut(G)}
\mbox{}\par
\begin{itemize}
\item[(i)]
$\ov{\ov A}=\N_{\Out(N)}(\ov{\ov H})$.
\item[(ii)]
$\Aut(G)\cong N.\N_{\Out(N)}(\ov{\ov H})$.
\end{itemize}
\end{lemma}
\proof (i)
By the above argument, we only need to show

\[\N_{\Out(N)}(\ov{\ov H})\leqslant\ov{\ov A}.\]

Let $\ov H=H\Z(N)/\Z(N)$.
Then $\ov G=\Inn(N){:}\ov H$.
The holomorph of $N$ is the semidirect product
$X:=R(N){:}\Aut(N)$ with $\Aut(N)$
acting naturally on $R(N)$.
The subgroup $R(N){:}\ov G$ of $X$
has a subgroup $\widehat{G}:=R(N){:}\ov H$,
and $\widehat{G}\cong\Inn(G)\cong G$ because
the action of $\ov H$ on $R(N)$ is by definition
the same as the action of $H$ on $N$.

For any $\ov\a\in\Out(N)$,
denote by $\a$ a preimage of $\ov\a$ in $\Aut(N)$,
and on the contrary for any $\b\in\Aut(N)$, denote
by $\ov\b$ the image of $\b$ in $\Out(N)$.

Let $\ov\theta\in\N_{\Out(N)}(\ov{\ov H})$.
Since $\ov{\ov H}^{\ov\theta}=\ov{\ov H}$
and $(\Inn(N))^\theta=\Inn(N)$,
we conclude that $\ov G^{\theta}=\ov G$.
Noting that $\ov H$ is a Hall subgroup of $\ov G$,
it follows from Schur-Zassenhaus's Theorem
that $\ov H^{\theta\nu}=\ov H$,
where $\nu\in\Inn(N)$.
This implies that
$\theta\nu$ normalises $\widehat{G}$
and so induces an automorphism
of $\widehat{G}$ by conjugation.
Thus $\theta\nu$ induces an automorphism $\vartheta$ of $G$
by an isomorphism from $G$ to $\widehat{G}$,
and hence $\vartheta\in A$.
It follows that $\theta\nu\in\ov A$
because $N^{\theta\nu}=N^{\vartheta}$.
Since $\nu\in\Inn(N)$, we have
$\ov{\theta\nu}=\ov\theta$, belonging to $\ov{\ov A}$.
Consequently, $\ov{\ov A}=\N_{\Out(N)}(\ov{\ov H})$, as required.

(ii) Since $\ov{\ov A}=\ov A/\Inn(N)\cong A/N$, we have
$A/N\cong\N_{\Out(N)}(\ov{\ov H})$ by part~(i).
It follows that $A\cong N.\N_{\Out(N)}(\ov{\ov H})$.
This completes the proof.
\qed

Finally, we quote a result about Maschke's theorem,
which will be used later.

\begin{lemma}\label{Mas}
{\rm (see~\cite[p.66]{Gorenstein})}
Let $V$ be a representation of the finite group $G$ over a
field $\FF$ in which $|G|$ is invertible.
Let $U$ be an invariant subspace of $V$. Then there exists
an invariant subspace $W$ of V such that $V=U\oplus W$ as representations.

\end{lemma}

\section{The automorphism groups}

In this section, we determine the automorphism groups
of Frobenius groups
\[\mbox{$G:=V{:}H=\ZZ_p^d{:}\ZZ_n$, where $p$ is a prime,
and $d,n\geqslant 1$}.\]

By Lemma~\ref{Aut(G)},
$\Aut(G)/V\cong\N_{\Aut(V)}(H)$.
Therefore, in order to determine $\Aut(G)$,
we only need to determine the normaliser
\[M:=\N_{\Aut(V)}(H)=\N_{\GL(d,p)}(H).\]

Now we regard $V$ as a faithful $\FF_pH$-module.
If $H$ is irreducible on $V$, then  the normaliser $M$ has
been determined via reference to \cite[Theorem~7.3, p.187]{B}:
\[M=\GammaL(1,p^d).\]
We thus assume that $H$ is reducible on $V$.
By Lemma~\ref{Mas}, $V$ can be decomposed as
\[V_H=U_1\oplus\cdots\oplus U_t\]
such that each $U_i$ is a faithful irreducible $\FF_p H$-module.
If the $U_i$ are pairwise isomorphic, 
then we call $V_H=U_1\oplus\cdots\oplus U_t$
a homogeneous decomposition of $V$.

In the case where $H$ is reducible and
homogeneous on $V$, the normaliser
$M$ is determined in the proof of \cite[Lemma~4.5]{JinWei}
(see the second paragraph of the proof on page 14).
For the completeness, we state it in the next lemma and give a short proof here.

By \cite{DF}, all the faithful irreducible $\FF_pH$-modules
have equal dimension, and moreover,
this dimension is the order of $p$ modulo $n$.
Throughout this article,
we use $f$ to denote
this same dimension. 

\begin{lemma}\label{homo}
Assume that $V$ is a homogeneous $\FF_pH$-module.
Then \[\mbox{$M=\GammaL(t,p^f)$ where $d=tf$}.\]
\end{lemma}
\proof
By our assumption, $V$ has the decomposition as
\[V_H=U_1\oplus\dots\oplus U_t\cong U^t,\]
where $U_i\cong U$ is a faithful irreducible $\FF_pH$-module.
Then $\dim_{\FF_p}U=f$ and  $tf=d$.
We can identify the action of $H$
on each $U_i$ with that of $H$ on $U$.
By \cite[27.14]{ASCH},
\[\C_{\GL(d,p)}(H)=\GL(t,p^f).\]
Since $\C_{\GL(d,p)}(H)\unlhd M$, it follows that
\[M\leqslant\N_{\GL(d,p)}(\C_{\GL(d,p)}(H))=\GL(t,p^f).\ZZ_f
=\GammaL(t,p^f),\]
where $\ZZ_f$ is the group of field automorphisms of $H$.
Let $Z\cong\ZZ_{p^f-1}$ be the center of $\GL(t,p^f)$. Then
$H\leqslant Z$, and thus $H$ \char\, $\GL(t,p^f)$.
It follows that $H\unlhd\GammaL(t,p^f)$.
Thus $M=\GammaL(t,p^f)$.
\qed

We now assume that the $\FF_pH$-module $V$ is not homogeneous.
Although $H$ is reducible, the overgroup $M=\N_{\GL(d,p)}(H)$ may be irreducible.
We proceed our proof in two subsections,
which treat the irreducible and reducible cases, respectively.

\subsection{The normaliser $\N_{\GL(d,p)}(H)$ is irreducible}\label{irreducible}

In this subsection, we assume that $M$ is irreducible on $V$.
By Clifford's Theorem,  $V$ can be decomposed as
\begin{equation}\label{equa-1}
V_H=W_1\oplus\cdots\oplus W_m,
\end{equation}
such that the $W_i$ are non-isomorphic homogeneous $\FF_pH$-modules. Since $H$ is normal in $M$, it follows that
$M$ preserves this
direct sum decomposition of $V$.
Thus
\[M\leqslant (\GL(W_1)\times\dots\times\GL(W_m)){:}\S_m,\]
where $\S_m=\Sym(\{1,\dots,m\})$.
By Lemma~\ref{homo}, all $W_i$ have the same
dimension over $\FF_p$, say $ef$.
In what follows, we will see that
$M$ can be embedded into $\GammaL(e,p^f)\wr\S_m$.

\vskip0.1in
Let $G_i$ be the quotient group of $G$ modulo $\prod_{j\not=i}W_j$
where $1\leqslant i\leqslant m$.
Then $G_i=W_i{:}H_i$ where $H_i\cong H=\ZZ_n$, and it implies that
$G_i=\ZZ_p^{ef}{:}\ZZ_n$ is a Frobenius group.
Let $H_i=\l h_i\r$ where $1\leqslant i\leqslant m$.
Then the group $G=(W_1\times\cdots\times W_m){:}H$
can be embedded  into
\begin{equation}\label{equa-2}
(W_1{:}H_1)\times\cdots\times (W_m{:}H_m),
\end{equation}
as a subgroup such that $H=\l h\r$
where $h=(h_1,\dots,h_m)$.

\vskip0.1in
Recall that $\dim_{\FF_p}W_i=ef$,
where $1\leqslant i\leqslant m$.
Since $\GL(W_i)\cong\GL(ef,p)$,
we may assume for convenience that $\GL(W_i)=\GL(ef,p)$.
Set 
\[B:=\GL(W_1)\times\cdots\times\GL(W_m)=\GL(ef,p)^m.\]

\begin{lemma}\label{LL}
$\C_{B\cap M}(H)=\GL(e,p^f)^m$ where $efm=\dim_{\FF_p}V$.
\end{lemma}
\proof By (\ref{equa-2}),
the action of $H_i$ on $W_i$ is equivalent to
the action of $H$ on $W_i$.
It follows that
$W_i$ is a homogeneous $\FF_pH_i$-module.
By Lemma~\ref{homo},
\[\mbox{$\C_{\GL(W_i)}(H_i)=\GL(e,p^f)$
where $1\leqslant i\leqslant m$}.\]

For any $c\in\C_{B\cap M}(H)$, we have
\[\mbox{$c=(c_1,\dots, c_m)$ where $c_i\in\GL(W_i)$}.\]
Then \[(h_1^{c_1},\dots, h_m^{c_m})=h^c=h=(h_1,\dots,h_m).\]
It follows that $h_i^{c_i}=h_i$, and so $c_i$ centralises $H_i$,
yielding $c_i\in\C_{\GL(W_i)}(H_i)$.
Thus $c\in\GL(e,p^f)^m$.  Consequently,
$\C_{B\cap M}(H)=\GL(e,p^f)^m$.
\qed

Let $L=L_1\times\cdots\times L_k$,
where $L_i\cong L_j$ for any $i,j$.
We call $\widehat L$ a diagonal subgroup of $L$ if
\[\widehat L=\{(x_1,x_1^{\varphi_2},\dots,x_1^{\varphi_k})\div x_1\in L_1\},\]
where $\varphi_i$ is an isomorphism from $L_1$ to $L_i$
for $2\leqslant i\leqslant k$.

By Lemma~\ref{LL}, $H$ can be embedded as a diagonal subgroup
into \[\GL(e,p^f)^m=\GL(e,p^f)\times\cdots\times\GL(e,p^f).\]
Since $B\cap M\unlhd M$, it follows that
$\C_{B\cap M}(H)\unlhd M$, and hence
\begin{equation}\label{equa-3}
\GL(e,p^f)^m\leqslant M\leqslant\GammaL(e,p^f)\wr\S_m.
\end{equation}

\begin{lemma}\label{Lei}
$B\cap M=\GL(e,p^f)^m.\ZZ_f$ where $efm=\dim_{\FF_p}V$.
\end{lemma}
\proof By (\ref{equa-2}), the action of
$H_i$ is equivalent to the action of $H$ on $W_i$.
It follows that
$W_i$ is a homogeneous $\FF_pH_i$-module.
By Lemma~\ref{homo},
\[\mbox{$M_i:=\N_{\GL(W_i)}(H_i)=
\GammaL(e,p^f)$ where $1\leqslant i\leqslant m$}.\]
We choose an element $\phi_i\in M$ such that
$h_i^{\phi_i}=h_i^p$. Then $M_i=\GL(e,p^f)\l\phi_i\r$.

For any $\s\in B\cap M$, we have
\[\mbox{$\s=(\s_1,\dots,\s_m)$ where $\s_i\in\GL(W_i)$}.\]
Since $\s$ normalises $H$,
there exists a positive integer $\lambda$ such that
\[(h_1^{\s_1},\dots, h_m^{\s_m})
=h^\s=h^\lambda=(h_1^\lambda,\dots, h_m^\lambda).\]
It follows that $h_i^{\s_i}=h_i^\lambda$,
and hence $\s_i\in M_i=\GammaL(e,p^f)$.
Thus \[\mbox{$\s_i=c_i\phi_i^{k_i}$,
where $c_i\in\GL(e,p^f)$, and $k_i$ is an integer}.\]
By the above equation, we have
\[h_i^\lambda=h_i^{\s_i}=h_i^{c_i\phi_i^{k_i}}=h_i^{p^{k_i}},\]
and hence
\[\mbox{$\lambda\equiv p^{k_i}\,(\mod n)$ where
$1\leqslant i\leqslant m$}.\]
It follows that
\[\mbox{$k_i\equiv k_j\equiv k\,(\mod f)$ for any $i,j$.}\]
Thus
$\s=(c_1\phi_1^k,\dots,c_m\phi_m^k)$.
So each element of $B\cap M$ can be written as
\[\mbox{$(c_1\phi_1^\ell,\dots,c_m\phi_m^\ell)$,
where $c_i\in\GL(e,p^f)$ and $\ell\geqslant0$}.\]
Let $\phi=(\phi_1,\dots,\phi_m)$.
Then $\phi\in B\cap M$. Therefore,
\[B\cap M=\GL(e,p^f)^m\l\phi\r=\GL(e,p^f)^m.\ZZ_f.\]
This completes the proof.
\qed

For convenience, we continue to use the
notation of Lemma~\ref{Lei}. By~(\ref{equa-3}),
elements of $M$ can be written as
\begin{equation}\label{equa-4}
\mbox{$c(\phi_1^{k_1},\dots,\phi_m^{k_m})\pi$,\,
where $c\in\GL(e,p^f)^m$,\, $k_i\geqslant0$\, and $\pi\in\S_m$}.
\end{equation}

\begin{lemma}\label{CEn}
$\C_{\GL(V)}(H)=\GL(e,p^f)^m$ where $efm=\dim_{\FF_p}V$.
\end{lemma}
\proof
Let $\s\in \C_{\GL(V)}(H)$. By~(\ref{equa-4}),
$\s$ can be written as
\[\mbox{$c(\phi_1^{k_1},\dots,\phi_m^{k_m})\pi^{-1}$,\,
where $c\in\GL(e,p^f)^m$,\, $k_i\geqslant0$\, and $\pi\in\S_m$.}\]
Suppose, if possible, that $\s\not\in\C_{B\cap M}(H)$.
If $\pi=1$, then there exists some $i_0$ for which $k_{i_0}\not\equiv0\,(\mod f)$,
where $1\leqslant i_0\leqslant m$.
We calculate that $h_{i_0}=h_{i_0}^{p^{k_{i_0}}}$,
which is absurd. Thus $\pi\not=1$.
So $\pi$ moves at least one point, say $r$.
Without loss of generality,
we may assume that $r^{\pi}=s$,
where $r\not=s$.
Calculation shows that
\[\begin{array}{rcl}
h^{\s}&=&(h_1^{\phi_1^{k_1}},\dots,
h_r^{\phi_r^{k_r}},\dots,
h_m^{\phi_m^{k_m}})^{\pi^{-1}}\\
&=& (h^{p^{k_{1^\pi}}}_{1^{\pi}},\dots,
h^{p^{k_{r^\pi}}}_{r^{\pi}},\dots,
h^{p^{k_{m^\pi}}}_{m^{\pi}})\\
&=&(h^{p^{k_{1^\pi}}}_{1^{\pi}},\dots,
h_s^{p^{k_s}},\dots,
h^{p^{k_{m^\pi}}}_{m^{\pi}}).
\end{array}\]
Noticing that $h^\s=h$,
we have $h_r=h_s^{p^{k_s}}$.

By (\ref{equa-2}), the action of $H_i$ on $W_i$ is
equivalent to the action of $H$ on $W_i$.
It follows that $W_i$ is a homogeneous $\FF_pH_i$-module.
Let $\widehat H=H_1\times\cdots\times H_m$.
Then $W_i$ is also a homogeneous $\FF_p\widehat H$-module.
Let $\hat h_i=(1,\dots,1,h_i,1,\dots,1)$ where
$1\leqslant i\leqslant m$.
Then
\begin{equation}\label{123}
h=(h_1,\dots,h_m)=\hat h_1\cdots\hat h_m.
\end{equation}
Let $(W_i)_{\widehat H}=U_i^e$, where $U_{i}\in\Irr(\FF_p\widehat H)$.
By Lemma~\ref{homo}, we may identify each $U_{i}$ with a field  $\FF{:}=\FF_{p^f}$ of order $p^f$,
and there exists $\o_i\in\FF^\times$ of order $n$
such that $\widehat H$ acts on each $x\in U_i$
by $\hat h_j\,{:}\,x\mapsto\o_i^{\d_{ij}}x$.
Since $h_r=h_s^{p^{k_s}}$,
we conclude that
$\o_r=\o_s^{p^{k_s}}$.

Let $\chi_i$ be the character defined by
$\FF_p\widehat H$-module $U_i$.
By \cite[25.10]{ASCH},
\[\chi_i(\hat h_i)=\o_i+\o_i^p+\cdots+\o_i^{p^{f-1}}.\]
By the previous paragraph, we have
$\chi_r(\hat h_r)=\chi_s(\hat h_s)$.
By (\ref{equa-2}) and (\ref{123}),
the action of $h$ on $W_i$ is the same as
the action of $\hat h_i$ on $W_i$.
It follows that $\chi_r(h)=\chi_s(h)$,
and hence $U_r\cong U_s$ ($\FF_p H$-isomorphic).
Thus $W_r$ and $W_s$ are isomorphic $\FF_p H$-modules.
This contradiction completes the proof.
\qed

\begin{proposition}\label{Nor}
$M=(\GL(e,p^f)^m.\ZZ_f).\ov L$,
where $\ZZ_f.\ov L\leqslant\Aut(H)$,
$\ZZ_f$ is the group of field automorphisms
of $H$, and $\ov L$ is abelian and regular on the set $\{1,2,\dots,m\}$.
\end{proposition}
\proof By (\ref{equa-1}), we have
\[V_H=W_1\oplus\cdots\oplus W_m,\]
where the $W_i$ are non-isomorphic homogeneous $\FF_pH$-modules.

Let $N=\C_{\GL(V)}(H).\ZZ_f$,
where $\ZZ_f$ is the group of field automorphisms
of $H$. By Lemma~\ref{Lei}, $N\unlhd M$.
For any element $\theta\in M\setminus N$,
we will show that
\[\mbox{$W_i^\theta\ncong W_i$ ($\FF_p H$-isomorphic)}.\]

Let $(W_i)_H=U_i^e$, where $U_i\in\Irr(\FF_p H)$.
Suppose that $W_i^\theta\cong W_i$ ($\FF_p H$-isomorphic).
Then $U_i^\theta\cong U_i$ ($\FF_p H$-isomorphic).
Let $\chi_i$ be the character defined by $\FF_pH$-module $U_i$.
Then $\chi_i^\theta$ is the character of $U_i^\theta$ over $\FF_p$.
Arguing as in Lemma~\ref{CEn}, \[\mbox{$\chi_i(h)=\sum_{j=1}^f\o^{p^{j-1}}$,
where $\o\in\FF_{p^f}^\times$ with $o(\o)=n$}.\]
Let $h^{\theta^{-1}}=h^\lambda$,
where $\gcd(\lambda,n)=1$. Then \[\chi_i^\theta(h)=\chi_i(h^{\theta^{-1}})=
\sum_{j=1}^f(\o^\lambda)^{p^{j-1}}.\]
Note that $\chi_i^\theta(h)=\chi_i(h)$.
By \cite[9.20]{Is}, we conclude that
$\lambda\equiv p^k\,(\mod n)$ for some $k$.
It follows that the action of $\theta$  on $H$
is induced by the field automorphism, and so
$\theta$ belongs to $N$, which is
a contradiction.
Thus $W_i^\theta\ncong W_i$,
as required.

Arguing similarly as above,
we may show that for $\vartheta\in N$,
$W_i^\vartheta\cong W_i$ ($\FF_p H$-isomorphic).
It follows that $N$ is the inertia
group of $W_i$ in $M$,
where $1\leqslant i\leqslant m$.

By Lemma~\ref{Lei},
$M=\C_{\GL(V)}(H).L$ and $\ZZ_f\leqslant L\leqslant\Aut(H)$.
Let $\ov L=M/N$ and
$[m]=\{1,2,\dots,m\}$.
By the previous paragraph, $M$
induces a permutation group $\ov L$ on $[m]$.
Since $V$ is a faithful irreducible $\FF_p M$-module,
we deduce that $\ov L$ is regular on $[m]$.
By Lemma~\ref{CEn},
$M=(\GL(e,p^f)^m.\ZZ_f).\ov L$,
competing the proof.
\qed

\noindent{\bf Remark.}
Although we have already characterized
the normaliser $\N_{\GL(V)}(H)$
(refer to Proposition~\ref{Nor}),
we can not completely
determine the structure of $\N_{\GL(V)}(H)$
because there are many different ways
to embed $H$ as a subgroup into $\GL(V)$ such that
$G$ is a Frobenius group when $|H|$ is arbitrarily large.

\vskip0.1in
We end this subsection by presenting several groups to explain the above remark.

\begin{example}
{\rm Let
$G=V{:}H\cong\ZZ_{31}^4{:}\ZZ_{15}$ be a Frobenius group.
By Lemma~\ref{Mas}, $V$ can be decomposed as
\[V=U_1\times U_2\times U_3\times U_4\]
such that $H$ normalises each $U_i$, and $U_i$ is
irreducible relative to the action of $H$.

Let $H=\l h\r$ with $h=(h_0^{r_1},h_0^{r_2},h_0^{r_3},h_0^{r_4})$,
where $o(h_0)=31$ and $\gcd(r_i,15)=1$.

\vskip0.07in
{\bf Case 1:}\, Suppose that $r_1=r_2=r_3=r_4$.

By Lemma~\ref{homo}, $\Aut(G)=\AGL(4,31)$.

\vskip0.07in
{\bf Case 2:}\, Suppose that $r_1=1$, $r_2=4$, $r_3=11$ and $r_4=14$.

Let $\pi_1=(12)(34)$ and $\pi_2=(13)(24)$.
Calculation shows that
\[\begin{array}{rcl}
h^{\pi_1}&=&(h_0^4,h_0,h_0^{14},h_0^{11})=h^4\,\,\, \mbox{and}\\
h^{\pi_2}&=&(h_0^{11},h_0^{14},h_0,h_0^4)=h^{-4}.
\end{array}\]
It follows that $\l\pi_1,\pi_2\r\leqslant\Aut(H)$.
By Proposition~\ref{Nor},
$\Aut(G)=\ZZ_{31}^4{:}(\ZZ_{30}^4{:}\ZZ_2^2)$.

\vskip0.07in
{\bf Case 3:}\, Suppose that $r_1=1$, $r_2=7$, $r_3=4$ and $r_4=13$.

Arguing as in Case~2, we have
$\Aut(G)=\ZZ_{31}^4{:}(\ZZ_{30}\wr\ZZ_4)$.

}
\end{example}

\subsection{The normaliser $\N_{\GL(d,p)}(H)$ is reducible}\label{3.4}
In this subsection, we begin with considering the case where
$V$ is a faithful reducible $\FF_p M$-module.

By Lemma~\ref{Mas}, $V$ can be decomposed as
\begin{equation}\label{equa-5}
V_H=U_1^{e_1}\oplus\cdots\oplus U_t^{e_t},
\end{equation}
where $U_i\in\Irr(\FF_pH)$, $U_i\ncong U_j$
and $e_i\geqslant1$.
As already mentioned above,
we continue to use  $f$ to denote
the same dimension of all the $U_i$.

\vskip0.07in
\noindent{\bf Remark.}
According to Proposition~\ref{Nor}, we obtain that $t\leqslant\varphi(n)/f$,
where $\varphi(n)$ is Euler's totient function.
\vskip0.07in

Let $G_{i}$ be the quotient group of $G$ modulo $\prod_{j\not=i}U_j^{e_j}$,
where $1\leqslant i\leqslant t$.
Then $G_{i}=U_i^{e_i}{:}H_i$ where $H_i\cong H=\ZZ_n$, and
it implies that $G_{i}=\ZZ_p^{e_if}{:}\ZZ_n$ is a Frobenius group.
Let $H_i=\l h_i\r$ where $1\leqslant i\leqslant t$.
Then the group $G=V{:}H=(U_1^{e_1}\times\cdots\times U_t^{e_t}){:}H$
can be embedded into
\begin{equation}\label{eauq1}
(U_1^{e_1}{:}H_1)\times\cdots\times(U_t^{e_t}{:}H_t),
\end{equation}
as a subgroup such that $H=\l h\r$ with $h=(h_1,\dots,h_t)$.

\begin{lemma}\label{center}
With the notation introduced above, the following holds.
\[\GL(e_1,p^f)\times\cdots\times\GL(e_t,p^f)\leqslant\C_{\GL(V)}(H).\]
\end{lemma}
\proof 
Let $\hat h_i=(1,\dots,1,h_i,1,\dots,1)$,
where $1\leqslant i\leqslant t$.
Then $h=\hat h_1\cdots\hat h_t$.
Let $\widehat U_i=U_i^{e_i}$.
By (\ref{eauq1}), we conclude that
the action of $\hat h_i$ on $\widehat U_i$ is the
same as the action of $h$ on $\widehat U_i$.
It follows that $\widehat U_i$ is
a homogeneous $\FF_pH_i$-module.
By Lemma~\ref{homo},
$\C_{\GL(\widehat U_i)}(H_i)=\GL(e_i,p^f)$.
Arguing as in Lemma~\ref{LL},
we can easily obtain the conclusion of this lemma.
\qed

Let $X\leqslant Y$, and let $V,U$ be $\FF X$-modules,
where $\FF$ is a field.
\begin{itemize}
\item[(i)]
Denote by $V^Y$ an induced module from $X$ to $Y$.

\vskip0.01in
\item[(ii)]
Write $U\div V$ if $V=V_1\oplus V_2$ with $U\cong V_1$,
where the $V_i$ are $\FF X$-submodules of $V$.

\vskip0.01in
\item[(iii)]
Denote by $n_V(U)$ the number of $W_i$
which are isomorphic to $U$ if $V=\oplus W_i$,
where the $W_i$ are $\FF X$-submodules of $V$.

\end{itemize}

\begin{lemma}\label{irrede}
$V_M=V_1\oplus\cdots\oplus V_s$, where
$V_i\in\Irr(\FF_pM)$ and $V_i\ncong V_j$ $(i\not=j)$.
\end{lemma}
\proof By (\ref{equa-5}), we have
\[\mbox{$V_H=U_1^{e_1}\oplus\cdots\oplus U_t^{e_t}$,
where $U_i\in\Irr(\FF_pH)$, $U_i\ncong U_j$ and $e_i\geqslant1$}.\]
Set \[\mbox{$T_i=
\{\theta\in M\div U_i^\theta\cong U_i\,
\mbox{($\FF_p H$-isomorphic)}\}$.}\]
By Mackey's Theorem (see \cite[p.174, Theorem~1.10]{Nagao}), $(U_i^{T_i})_H\cong U_i^{f_i}$
where $f_i\leqslant e_i$.
By \cite[12.12]{ASCH}, $\C_{\GL(V)}(H)$ is a subgroup of $T_i$.
Thus, by Lemma~\ref{center},
$\GL(e_i,p^f)\leqslant T_i$.
It follows that
$\Hom_{\FF_pH}(U_i,(U_i^{T_i})_H)\cong\FF_{p^f}^{e_i}$,
and hence $f_i=e_i$.
It further implies that $U_i^{T_i}$ is an irreducible $\FF_pT_i$-module.
By \cite[6.11]{Is}, $U_i^M$ is an irreducible $\FF_pM$-module.
As a convenience, we denote  by $V_i$ that module.
It is easily shown that $V_i$ is
isomorphic to an irreducible submodule of $V_M$.
Thus, without loss of generality,
we may identify $V_i$ with a submodule of $V_M$.

Next we show that
\[\mbox{if $U_j\div V_i$, then $n_{V_i}(U_j)=e_j=e_i$}.\]
By Clifford's Theorem, we conclude that
$n_{V_i}(U_j)=e_i$.
Arguing similarly as the first paragraph, there exists an
irreducible $\FF_pM$-module $V_j$
for which $V_j=U_j^M$ and
$\Hom_{\FF_pH}(U_j,(V_j)_H)\cong\FF_{p^f}^{e_j}$.
It follows that $\Hom_{\FF_pM}(V_j,V_i)\not=0$,
and so $V_j\cong V_i$ ($\FF_p M$-isomorphic).
Thus $n_{V_i}(U_j)=e_i=e_j$, as desired.

By our assumption,
$V$ is a faithful reducible $\FF_p M$-module.
Repeating the above process, we obtain
all the irreducible $\FF_p M$-submodules of $V$.
Without loss of generality, we may assume that $V_1,\dots, V_s$ are
all pairwise non-isomorphic irreducible $\FF_p M$-submodules of $V$.
Let $W=V_1\oplus\dots\oplus V_s$.
Suppose, if possible, that
$W$ is a proper submodule of $V$.
By the second paragraph, we may assume that
\[\mbox{$W_H=U_1^{e_1}\oplus\cdots\oplus U_{k}^{e_{k}}$, where $k<t$.}\]
Arguing as the first paragraph,
there exists an irreducible
$\FF_p M$-submodule of $V_M$, say $U$
such that $U_t^M=U$.
Since $\Hom_{\FF_pH}(U_t,(V_i)_H)=0$, it follows from
Frobenius reciprocity that $\Hom_{\FF_pM}(U,V_i)=0$ where $1\leqslant i\leqslant s$,
contrary to our assumption.
Thus $V_M=V_1\oplus\cdots\oplus V_s$.
This completes the proof.
\qed

Let $G_{i}$ be the factor group of $G$ modulo $\prod_{j\not=i}V_j$, where $1\leqslant i\leqslant s$.
Then $G_{i}=V_i{:}H_i$ where $H_i\cong H=\ZZ_n$,
and it implies that $G_{i}$
is a Frobenius group. Let $H_i=\l h_i\r$
where $1\leqslant i\leqslant s$.
Then the group $G=(V_1\times\cdots\times V_s){:}H$
can be embedded into
\[(V_1{:}H_1)\times\cdots\times(V_s{:}H_s),\]
as a subgroup such that $H=\l h\r$ with $h=(h_1,\dots,h_s)$.

\begin{lemma}\label{WW}
$M\leqslant\N_{\GL(V_1)}(H_1)\times\cdots\times\N_{\GL(V_s)}(H_s)$.
\end{lemma}
\proof By Lemma~\ref{irrede},  $V_M=V_1\oplus\cdots\oplus V_s$,
where $V_i\in\Irr(\FF_pM)$ and $V_i\ncong V_j$ for $i\not=j$.
Since $M$ fixes each $V_i$,
it follows that
\[M\leqslant\GL(V_1)\times\cdots\times\GL(V_s).\]
For any $\s\in M$, we may write $\s=(\s_1,\dots,\s_s)$
where $\s_i\in\GL(V_i)$.
Then there exists a positive integer $\lambda$
such that $h^\s=h^\lambda$, namely,
\[(h_1^{\s_1},\dots, h_s^{\s_s})=h^\s=h^\lambda
=(h_1^\lambda,\dots, h_s^\lambda).\]
It follows that $h_i^{\s_i}=h_i^\lambda$, and so
$\s_i$ normalises $H_i$, forcing that
$\s_i\in\N_{\GL(V_i)}(H_i)$.
Thus $\s\in\prod_{i=1}^s\N_{\GL(V_i)}(H_i)$.
Therefore, the lemma follows.
\qed

\begin{lemma}\label{WWWW}
$\C_{\GL(V)}(H)=\C_{\GL(V_1)}(H_1)\times\cdots\times\C_{\GL(V_s)}(H_s)$.
\end{lemma}
\proof
Take any $c\in\C_{\GL(V)}(H)$. By Lemma~\ref{WW},
\[c=(c_1,\dots,c_s)\,\,\mbox{where}\,\,c_i\in\GL(V_i).\]
Then we have
\[(h_1^{c_1},\dots, h_s^{c_s})=h^c=h=(h_1,\dots, h_s).\]
It follows that $h_i^{c_i}=h_i$, and so
$c_i$ centralises $H_i$, forcing that $c_i\in\C_{\GL(V_i)}(H_i)$.
Thus $c\in\prod_{i=1}^s\C_{\GL(V_i)}(H_i)$.
Since $c$ is arbitrary, we have
\[\C_{\GL(V)}(H)\leqslant\C_{\GL(V_1)}(H_1)
\times\cdots\times\C_{\GL(V_s)}(H_s).\]
It is clear that $\prod_{i=1}^s\C_{\GL(V_i)}(H_i)\leqslant\C_{\GL(V)}(H)$.
Thus the lemma holds.
\qed

Applying Lemmas~\ref{CEn} and~\ref{WWWW},
we have the following corollary.

\begin{corollary}\label{Cen}
With the hypothesis of Subsection~$3.2$, we have
\[\C_{\GL(V)}(H)=\GL(e_1,p^f)\times\cdots\times\GL(e_t,p^f).\]
\end{corollary}

Let $\rho_i$ ($1\leqslant i\leqslant s$) be the projection map
\begin{equation}\label{equa-7}
\prod_{j=1}^s\N_{\GL(V_j)}(H_j)\to \N_{\GL(V_i)}(H_i):(x_1,\dots,x_s)\mapsto(1,\dots,1,x_i,1,\dots,1).
\end{equation}
By Lemma~\ref{WW}, we use $M_i$ for the image
$M^{\rho_i}$ where $1\leqslant i\leqslant s$.

\begin{lemma}\label{WWW}
With the notation above,
$M_i=\N_{\GL(V_i)}(H_i)$ where $1\leqslant i\leqslant s$.
\end{lemma}
\proof 
Let $N_i=\N_{\GL(V_i)}(H_i)$ where $1\leqslant i\leqslant s$.
Then $M_i\leqslant N_i$.
By (\ref{equa-7}),
the action of $M_i$ on $V_i$ is equivalent to
the action of $M$ on $V_i$.
It follows that $M_i$ acts
irreducibly on $V_i$.
So does $N_i$.
By Proposition~\ref{Nor},  $N_i=\C_{\GL(V_i)}(H_i).L_i$, and
$\l\hat\phi_i\r\leqslant L_i\leqslant\Aut(H_i)$, where
$\hat\phi_i=(1,\dots,1,\phi_i,1,\dots,1)$ is a field automorphism of $H_i$ of order $f$.
Without loss of generality, we may assume
$h_i^{\phi_i}=h_i^p$, where $1\leqslant i\leqslant s$.

For convenience, we again denote by $\hat\phi_i$ the
preimage of $\hat\phi_i$ in $N_i$ under $N_i\to L_i$.
Set $\phi=(\phi_1,\dots,\phi_s)$. Then
\[h^\phi=(h_1^{\phi_1},\dots,h_s^{\phi_s})=(h_1^p,\dots,h_s^p)=h^p.\]
It follows that $\phi\in M$.
By (\ref{equa-7}), we obtain $\hat\phi_i\in M_i$ .
Let $K_i=\C_{\GL(V_i)}(H_i)\l\hat\phi_i\r$.
By Lemma~\ref{WWWW},  $K_i\leqslant M_i$.
By Clifford's Theorem,
\[(V_i)_{H_i}=W_{i1}\oplus\cdots\oplus W_{i\ell_i},\]
where the $W_{ij}$ are non-isomorphic
homogeneous $\FF_pH_i$-modules.
By Proposition~\ref{Nor}, $K_i$ is
the inertia group of $W_{ij}$ in $N_i$, and $\ell_i=|L_i|/f$.
Since $M_i$ is irreducible on $V_i$,
we conclude that $M_i$ induces a transitive permutation group
$M_i/K_i$ on $\{W_{i1},\dots,W_{i\ell_i}\}$.
It follows that $|M_i/K_i|=\ell_i$.
Thus $M_i=N_i$ where $1\leqslant i\leqslant s$.
This completes the proof of Lemma~\ref{WWW}.
\qed

\begin{lemma}\label{FM}
Let $V_i$ be as in Lemma~$\ref{irrede}$
where $1\leqslant i\leqslant s$.
Then
\begin{itemize}
\item[(i)]
$(V_i)_{H}=V_{i1}^{e_i}\oplus\cdots\oplus V_{i\ell}^{e_i}$ and
$\dim_{\FF_p}V_{ij}=f$,
where $V_{ij}\in\Irr(\FF_pH)$, $V_{ij}\not\cong V_{ik}$ $(j\not=k)$, and $f=\ord_n(p)$;
\item[(ii)]
$M=((\GL(e_1,p^f)^\ell\times\cdots\times\GL(e_s,p^f)^\ell).\ZZ_f).L$,
where $\ZZ_f.L\leqslant\Aut(H)$, $|L|=\ell$,
and $(e_1+\dots+e_s)f\ell=d$.
\end{itemize}
\end{lemma}
\proof 
Since $H\unlhd M$, it follows from Clifford's Theorem that
\[(V_i)_H=V_{i1}^{e_i}\oplus\cdots\oplus V_{i\ell_i}^{e_i},\]
where 
$V_{ij}\in\Irr(\FF_pH)$, $V_{ij}\ncong V_{ik}$ ($j\not=k$)
and $\ell_i\geqslant1$.
By \cite{DF}, $\dim_{\FF_p}V_{ij}=f$ for all $i,j$.

Recall that $\rho_i$ is the projection map defined
in the paragraph after Corollary~\ref{Cen},
where $1\leqslant i\leqslant s$.
Then $M_i=M^{\rho_i}$ and $H_i=H^{\rho_i}$.
By (\ref{equa-7}), the action of $M_i$ on $V_i$ is
equivalent to the action of $M$ on $V_i$.
Therefore, we may identify the $V_i$ with
an irreducible $\FF_pM_i$-module.
Similarly, we may also identify the $V_{ij}$
with non-isomorphic faithful irreducible $\FF_pH_i$-modules.
By Lemma~\ref{WWW} along with Proposition~\ref{Nor},
\[\mbox{$M_i=\C_{\GL(V_i)}(H_i).L_i$\,
and $\l\hat\phi_i\r\leqslant L_i\leqslant\Aut(H_i)$},\]
where $\hat\phi_i=(1,\dots,1,\phi_i,1,\dots,1)$ is a field automorphism of $H_i$ of order $f$, and
$|L_i|=f\ell_i$.
Let $C=\C_{\GL(V)}(H)$ and $C_i=\C_{\GL(V_i)}(H_i)$.
By Lemma~\ref{WWWW}, $C^{\rho_i}=C_i$,
and by definition,
$L_i=M_i/C_i=M^{\rho_i}/C^{\rho_i}$,
where $1\leqslant i\leqslant s$.
Let $\ov M=M/\C_{\GL(V)}(H)$.
Noting that $M\leqslant M_1\times\cdots\times M_s$,
it follows that
$\ov M\lesssim L_1\times\cdots\times L_s$.
Since $\ov M\leqslant\Aut(H)$, we deduce that
$\ov M$ is isomorphic to a diagonal subgroup of $\prod_{i=1}^sL_i$.
It follows that
$\ov M\cong L_i\cong L_j$, and
hence $\ell_i=\ell_j=\ell$,
where $1\leqslant i,j\leqslant s$,
as in Lemma~\ref{FM}~(i).

Let $\phi=(\phi_1,\dots,\phi_s)$ be as in Lemma~\ref{WWW}.
Then $\phi$ induces an automorphism of $H$ of order $f$.
Let $N=\C_{\GL(V)}(H)\l\phi\r$.
Then $N\unlhd M$.
Let $\ov N=N/\C_{\GL(V)}(H)$.
Then $|\ov N|=f$.
Since $M/N\cong\ov M/\ov N$ and $|\ov M|=f\ell$,
it follows that $|M/N|=\ell$.
By Proposition~\ref{CEn},
$\C_{\GL(V_i)}(H_i)=\GL(e_i,p^f)^{\ell}$
where $1\leqslant i\leqslant s$.
By Corollary~\ref{Cen},
\[\C_{\GL(V)}(H)=\prod_{i=1}^s\C_{\GL(V_i)}(H_i)
=\prod_{i=1}^s\GL(e_i,p^f)^\ell.\]
Let $M/N=L$. By the previous argument,
$M=((\prod_{i=1}^s\GL(e_i,p^f)^\ell).\ZZ_f).L$,
where $(\sum_{i=1}^s e_i)f\ell=d$,
as in Lemma~\ref{FM}~(ii),
completing the proof.
\qed

\vskip0.1in
The assertion of Theorem~\ref{WLLLL} follows from
Lemma~\ref{Aut(G)} and Lemma~\ref{FM}.

\vskip0.1in
In particular, Theorem~\ref{WLLLL} implies the following consequences.

\begin{corollary}\label{Final}
Keep the notation and conditions of Theorem~$\ref{WLLLL}$. Then
\[\C_{\GL(d,p)}(H)=\GL(e_1,p^f)^\ell\times\cdots\times \GL(e_s,p^f)^\ell.\]
\end{corollary}

\vskip0.07in
For $X\leqslant Y$ define $\Aut_Y(X)=\N_Y(X)/\C_Y(X)$ to be the automizer in $Y$ of $X$.
Then $\Aut_Y(X)\leqslant\Aut(X)$, and indeed $\Aut_Y(X)$ is the group of automorphisms
induced on $X$ in $Y$.

\begin{corollary}\label{c}
Let $G=V{:}H=\ZZ_p^d{:}\ZZ_{q^\ell}$ be a Frobenius group, where $p,q$ are two distinct primes,
and $d,\ell$ are positive integers.
Let $f=\ord_{q^\ell}(p)$. Then the following hold.
\begin{itemize}
\item[(i)]
If $q=2$, then $\Aut_{\GL(d,p)}(H)=L$
where $\ZZ_f\leqslant L\leqslant\ZZ_2\times\ZZ_{2^{\ell-2}}$.

\item[(ii)]
If $q$ is odd, then $\Aut_{\GL(d,p)}(H)=\ZZ_k$
where $\ZZ_f\leqslant\ZZ_k\leqslant\ZZ_{q^{\ell-1}(q-1)}$.
\end{itemize}
\end{corollary}

\vskip0.2in
{\noindent\bf Acknowledgements.}
This paper is a part of the author's
PhD thesis under the
supervision of Professor Caiheng Li.
The author would like to thank Professor
Caiheng Li for his stimulating discussions and many
helpful suggestions.

\end{document}